\documentclass[12pt]{article}

\usepackage{amsfonts}
\usepackage{amssymb}
\usepackage{amsmath}
\usepackage[T1]{fontenc}
\usepackage{graphicx}
\usepackage{latexsym}

\def\cG{{\mathcal G}}

\def\cS{{\mathcal S}}
\def\cC{{\mathcal C}}

\def\pl{\partial}

\def\Llrw{\Longleftrightarrow}
\def\wt{\widetilde}
\def\bk{\backslash}

\newlength{\boxwidth}
\def\UN{\hbox{ 1 \hskip -7pt I}}

\def\N{{\mathbb{N}}}
\def\R{{\mathbb R}}
\def\Z{{\mathbb Z}}

\def\P{{\mathbb P}}
\def\E{{\mathbb E}}
\def\T{{\mathbb T}}
\def\]{{\mathbb ]}}
\def\[{{\mathbb [}}

\def\exp{\mathrm{exp}}
\def\e{\mathrm{e}}

\newcommand{\fdefeq}{\overset{\text{def.}}{=}}

\newtheorem{thm}{Theorem}[section]     
\newtheorem{pro}[thm]{Proposition}

\newtheorem{lem}[thm]{Lemma}
\newenvironment{proof} {\noindent {\em \bf{Proof}} \\ }{ \hfill $\Box$  }

\newtheorem{remk}[thm]{Remark}     

\def\build#1_#2^#3{\mathrel{\mathop{\kern 0pt#1}\limits_{#2}^{#3}}}

\title{Internal Diffusion Limited Aggregation on discrete groups having exponential growth}

\author{S\'ebastien Blach\`ere\footnote{Research supported by the ESF program 'Phase transition and fluctuation phenomena for random dynamics in spatially extended systems' ref. 554 and by Marie Curie Fellowship HPMF-CT-2002-02137}, Sara Brofferio\footnote{Research supported by  Marie Curie Fellowship HPMF-CT-2002-02137}}

\date{}

\begin{document}

\maketitle

\footnotetext{{\bf Keywords}: Interacting particle systems, random walks on groups, Green function}

\footnotetext{{\bf AMS 2000 subject classification}: 60B15, 60K35, 82B24}

\unitlength=1cm

\begin{abstract}
  The Internal Diffusion Limited Aggregation has been introduced by
Diaconis and Fulton in 1991. It is a growth model defined on an infinite
set and associated to a Markov chain on this set. We focus here on sets which are finitely generated groups with exponential growth.
We prove a shape theorem for the Internal DLA on such groups associated to symmetric random walks.
For that purpose, we introduce a new distance associated to the Green function, which happens to have some interesting properties.
In the case of homogeneous trees, we also get the right order for the fluctuations of that model around its limiting shape.
\end{abstract}

\section{Introduction}


Let $\cG$ be an infinite discrete set and let $P$ be a irreducible transition matrix on that set.
The matrix $P$ allows to define a sequence of i.i.d. Markov chains $(S^n (\cdot ))_{n\in \N^*}$ with a common starting point $O$ ($S(0)=O$).
The Internal Diffusion Limited Aggregation (Internal DLA) is a Markov chain $A(n)$ $(n\in \N)$ of increasing subsets of $\cG$ defined as follows:
\begin{itemize}
\item[$\bullet$] $A(1)= \{ O\}$,
\item[$\bullet$] $A(n+1)=A(n)\cup S^{n+1}(\tau_{A(n)^c})$,
\end{itemize}
where $\tau_A$ is the hitting time of the set $A$.

Thus, the set $A(n)$ grows as follows.
At each discrete time $n$, we start the Markov chain $S^n(\cdot )$ (moving according to $P$) at $O$, wait until it leaves the previous set $A(n-1)$, and add the first element
visited outside $A(n-1)$, to obtain $A(n)$.
Note that the Markov chain $A(n)$ is well defined since, for all $n$, $\tau_{A(n)^c}<\infty$ almost surely, as $A(n)$ is finite and $S^{n+1}(\cdot )$ is irreducible.

The name Internal DLA comes from the similitude with the (external) DLA. In the latter model, the Markov chains start ``at infinity'', conditioned to hit the previous cluster, and are stuck just before hitting it. See for example Lawler \cite{lawler91} for a precise definition and some properties.

The Internal DLA was introduced by Diaconis and Fulton in 1991 \cite{diaconis91}. They gave the first shape theorem on $\Z$ where $S^n(\cdot )$ are Simple Random Walks (nearest neighbours random walks with uniform transition probability).
Later, other shape theorems were proved on $\Z^d$ (Lawler, Bramson \& Griffeath \cite{lawler91}, Lawler \cite{lawler95}, Blach\`ere \cite{seb02}).
These results set that, when the random walk is centered (under some moment condition), the model has a limiting shape (see one of these for precise definition).
Roughly, that means that the shape of $A(n)$  tends to the shape of the balls of $\R^d$ for a norm associated to the transition matrix $P$.


In the present article, we focus on groups $\Gamma$ having exponential growth.
It means that the balls (associated to the word metric for any finite generating set) have cardinality growing exponentially in the radius.
On such groups, we consider symmetric random walks with finitely supported increment.
Our main result (Theorem \ref{limshape}) is a shape theorem for the Internal DLA.
The meaning of a shape theorem is to compare the random shape of the Internal DLA with a deterministic exhaustion of the group.
For that purpose, we define a new distance on groups, called the hitting distance, by
$$
d(x,y)\fdefeq - \ln (F(x,y)) \, ,
$$
where $F(x,y)$ is the hitting probability of $y$ for a random walk starting at $x$.
As the Green function $G(\cdot , \cdot )$ is proportional to the function $F(\cdot , \cdot)$, the balls for that distance are in fact the sets of the form,
$$
\{ x \; \mathrm{s.t.} \; G(e,x)\geq N\} \, ,
$$
where $e$ is the identity of the group.
These sets exhaust $\Gamma$ when $N$ goes to $0$.

This distance can be compared to the word distance but they are
not equivalent in general. Outside its link with the Internal DLA,
an interesting property of that distance is that one can deduce
the exact  exponential growth rate of the balls (Proposition
\ref{voln}).

The shape theorem states that the Internal DLA is close to the
balls $B(n)$ of radius $Kn$ for the hitting distance ($K$ is a
constant that ensures that the ball $B(n)$ contains the boundary
of the ball $B(n-1)$). More precisely there exists some constants
$C_I$ and $C_O$ such that at time $V(n)$ (=volume of $B(n)$)  the
sets of the IDLA satisfy:
$$
B(n-C_I \ln n) \subset A(V(n)) \subset B (n+C_O \sqrt{n} )
$$
for every large $n$, with probability $1$.
 This result gives the
upper bounds for the fluctuations between $A(V(n))$ and $B(n)$:
the inner fluctuations are at most logarithmic, whence the outer
fluctuations are at most of order square root of the radius. The
proof of that result is adapted from the one on $\Z^d$ (Lawler,
Bramson \& Griffeath \cite{lawler92} and Lawler \cite{lawler95})
with new ingredients related to the hitting distance and some
simplifications, partially due to the setting of groups with
exponential growth.

For the Simple Random Walk of homogeneous trees (where the hitting
distance and the word distance are proportional) such a result was
already obtained by Etheridge and Lawler (unpublished). In that
setting, they also announced (in \cite{lawler91}) that the lower bounds for the inner
and outer fluctuations were of the same order than the upper ones,
but according to Lawler, a written proof of that result was
missing and the lower bound for the inner error could be
considered  an open question. Therefore, we prove (Propositions
\ref{lboutTq} and \ref{lbinnTq}) that in the case of a homogeneous
tree  the order of the inner fluctuation cannot be smaller of $\ln
n$, while the outer fluctuation is bigger than some constant times
$\sqrt n$.

The article is constructed as follows. Section \ref{symRW} defines
the setting of symmetric random  walks on groups with exponential
growth, defines the ``hitting distance'' and compares it to the
word distance. Section \ref{shapeth} gives the proof of our main
result, the shape theorem for the Internal DLA with upper bounds
for the fluctuations. Section \ref{trees} focus on Simple Random
Walk on homogeneous trees and proves the sharpness of the bounds
for these examples.

\medskip

We would like to warmly acknowledge Wolfgang Woess for his help and support (both mathematical and financial) along this work.


\section{Symmetric random walk on discrete groups with exponential growth} \label{symRW}

Let $\Gamma$ be a finitely generated group. Let $\cS$ be a finite
symmetric set of generators. We write each element $x$ of $\Gamma$
as a word composed by these generators (letters). The minimal
length among the words representing $x$ is called the word
distance $|x|$. Let $e$ be the identity of $\Gamma$. According to
this distance, we define the balls $B_w(n)=B_w(e,n)\fdefeq \{ x\in
\Gamma \; : \;  |x| \leq n \}$ and their volume $V_w(n)\fdefeq \#
B_w(n)$. The exponential growth property corresponds to
$$
V_w(n) \geq \exp (c\, n) \, ,
$$
for some constant $c$. Note that an exponential upper bound of the
volume is always valid since the fastest  volume growth occurs
when there is no relation between the generators. Hence, $V_w (n)
\leq 1 + (\# \cS -1 )^n$. Note also that this property does not
depend on the choice of the generating set.

Let $\mu$ be a symmetric probability law whose support is $\cS$.
Let $(X_k)_{k\in \N^*}$ be a sequence of i.i.d. random variables
whose common law is $\mu$. The process $$S(k)\fdefeq xX_1X_2\cdots
X_k,$$ with $S(0)= x$, is a symmetric irreducible random walk on
$\Gamma$  starting at $x$. The process $S(\cdot )$ is then a
nearest neighbors symmetric random walk on the Cayley graph
$\cG=(\Gamma , \cS)$ on which we study the Internal DLA (with
$A(0)=\{ e\}$). We denote $\P^x$ and $\E^x$, respectively, the
probability and expectation related to a random walk starting at
$x$. When $x=e$, the exponent will be omitted.

The Green function $G(x,y)$ is defined as the expected number of visits at $y$ for a random walk starting at $x$:
$$
G(x,y)\fdefeq  \E^x \left[ \sum_{k=0}^{\infty} \UN_{\{S(k)=y\}} \right] = \sum_{k=0}^{\infty} \P^x [ S(k)=y ] \, .
$$
Since the group has an exponential growth,  every random walk is
transient and so the above summation is finite ($G(x,y)$ is
actually bounded uniformly in $x$ and $y$).

For a random walk $S(\cdot )$, let $\tau_y$ be the first hitting time of $y$:
$$
\tau_y \fdefeq \inf \{ k\geq 0 \; : \; S(k) = y\} \, .
$$
When $y$ is never attained, let $\tau_y=\infty$.
The hitting probability of $y$ starting at $x$ is denoted
$$
F(x,y) \fdefeq  \P^x [ \tau_y < \infty ] \, .
$$
We write $G(x)=G(e,x)$ and $F(x)=F(e,x)$.
Note that $F$ and $G$ are symmetric and invariant by left multiplication.
In particular, $G(y,y)=G(e)$.
Therefore, the functions $F$ and $G$ are proportional since a straightforward computation shows
\begin{equation} \label{GF}
G(x,y) = G(y,y) F(x,y) = G(e) F(x,y) \, .
\end{equation}


We are now ready to define the distance we will use to study the Internal DLA.
For all $x,y \in \Gamma$, we define
$$
d(x,y) \fdefeq - \ln F(x,y) \, .
$$
Since $F(x,y)$ is the hitting probability of $y$ starting at $x$,
we call $d(\cdot , \cdot)$  the \textbf{hitting distance}.

\begin{lem} \label{dist}
$d(\cdot , \cdot )$ is a left invariant distance on $\Gamma$.
\end{lem}

\begin{proof}
As $F(x,y) \leq 1$, $d(\cdot , \cdot)$ is always non-negative.

The symmetry of the random walk and the left invariance of $F(\cdot , \cdot)$ yield the same properties for $d(\cdot , \cdot)$.

Remark also that since the random walk is symmetric and transient
$$
\forall x\neq y , \; 1 > \P^x [\tau_x' <\infty] \geq \P^x [\tau_y <\infty] \P^y [\tau_x <\infty] =F(x,y)^2\, ,
$$
where $\tau_x' \fdefeq \inf \{ k\geq 1 \; : \; S(k) =x\}$. Thus
$$
d(x,y)=0 \Llrw F(x,y)=1 \Llrw x=y \, .
$$
Finally,
$$
\P^x [\tau_z <\infty] \geq \P^x [\tau_y <\infty] \P^y [\tau_z <\infty]
$$
leads to the triangular inequality.
\end{proof}

Let define
\begin{equation} \label{defK}
K\fdefeq \max \{ d(e,x) \; : \; x \in \cS\} \, .
\end{equation}
According to this distance, we define the ball of radius $Kn$,
$$
B(x,n)\fdefeq  \{ y \in \Gamma \; : \; d(x,y) \leq Kn \}
$$
and its external boundary
$$
\partial B(x,n) \fdefeq \{ z \not\in B(x,n) \; : \; \exists y\sim z \mbox{ and } y \in B(x,n) \} \, .
$$
Here $y\sim z$ means that $y$ and $z$ are neighbours, that is $y^{-1}z \in {\cS}$.

The reason why the constant $K$ appears in the definition of the balls is to be sure that $\partial B(e,n) \subset B(e,n+1)$.
Indeed, let $x\in \partial B(e,n)$ and $y \in B(e,n)$ such that $y \sim x$.
Then,
$$
d(e,x)\leq d(e,y)+d(y,x) \leq Kn + \max \{ d(e,z) \; : \; z \in \cS\} \leq K(n+1) \, .
$$
We write $d(x)=d(e,x)$ and $B(n)=B(e,n)$.
To have a better understanding of this distance, we will compare it to the word distance.
First remark that, for any $x\in \Gamma$, writing $x=x_1\cdots x_{|x|}$ with all the $x_i$'s in $\cS$, leads to
\begin{equation} \label{dGleqdw}
d(x) \leq \sum_{i=1}^{|x|} d(x_i) \leq K |x| \, .
\end{equation}
To get a lower bound for the hitting distance, we need asymptotic estimates of the Green function.
These estimates are direct consequences of the following results on the asymptotic of the iterated transition probability:

\noindent 1. An off-diagonal gaussian upper bound for symmetric
random walks (with finite support) on generic groups (Varopoulos
\cite{varo85a}; Carne \cite{carne85}),
\begin{equation} \label{gaussianP}
\exists C_g>1 \mbox{ s.t. } \forall x,y \in \Gamma \,\mbox{ and } k\in\N  \;\;\; \P^x [S(k) =y ] \leq C_g \exp (- C_g |x^{-1}y|^2/k ) \, ,
\end{equation}
\noindent 2. For random walks on groups having exponential growth (Varopoulos \cite{varo85b}; Carlen, Kusuoka \& Stroock \cite{carlen87}),
\begin{equation} \label{returnP}
\exists C_{e}>1 \mbox{ s.t. } \forall x,y \in \Gamma \, \mbox{ and } k\in\N\;\;\; \P^x [S(k) =y ] \leq C_{e} \exp (- C_{e} k^{1/3} ) \, ,
\end{equation}
\noindent 3. In the special case of non-ameanable groups, the latter estimates becomes (Kesten \cite{kesten59})
\begin{equation} \label{returnPna}
\exists C_{na}>1 \mbox{ s.t. } \forall x,y \in \Gamma \, \mbox{ and } k\in\N\;\;\; \P^x [S(k) =y ] \leq C_{na} \exp (- C_{na} k ) \, .
\end{equation}
For an overview of such estimates and precise references, see
Woess \cite{woess00} or Pittet \& Saloff-Coste \cite{pittet}.

\begin{lem} \label{BsupGw}
Let $\Gamma$ be a group having exponential growth, then $\exists
C>0$ such that
$$
G(x) \leq \exp ( -C|x|^{1/2} ) \;\; \forall x\in \Gamma.
$$
\end{lem}

\begin{proof}
Using (\ref{gaussianP}) and (\ref{returnP}),
\begin{eqnarray*}
G(x) & = & \sum_{k=0}^{\infty} \P [ S(k) =x] =  \sum_{k \leq |x|^{3/2}} \P [ S(k) =x] + \sum_{k>|x|^{3/2}} \P [ S(k) =x]\\
 & \leq & \sum_{k \leq |x|^{3/2}} C_g \exp (- C_g |x|^2/k ) + \sum_{k> |x|^{3/2}} C_e \exp (- C_e k^{1/3} ) \\
 & \leq & C_g |x|^{3/2} \exp ( -C_g |x|^{1/2} ) + 15  |x| \exp ( -C_e |x|^{1/2} )\\
 & \leq & \exp ( -C|x|^{1/2} ) \, ,
\end{eqnarray*}
for some constant $C$. The before last inequality comes from the
fact that for any sufficiently large $N$
\begin{eqnarray*}
    \sum_{k=N+1}^\infty \exp(- C_e k^{1/3})
    &\leq & \int_N^{+\infty}\exp(- C_e t^{1/3}) dt\\
    & \leq & \int_N^{+\infty}\frac{d(-6 t^{2/3}\exp(- C_e
    t^{1/3}))}{dt}dt= -6 N^{2/3}\exp(- C_e
    N^{1/3})
\end{eqnarray*}
\end{proof}

Actually, Revelle \cite{revelle05} proved independently Lemma \ref{BsupGw}.

Note that when the group is non-ameanable, using (\ref{returnPna}), this upper bound of the Green function becomes:
\begin{equation} \label{BsupGwna}
G(x) \leq \exp ( -C|x|) \, .
\end{equation}

The above lemma assures that $G(x)$ tends to $0$ as $|x|$ goes to infinity.
By (\ref{GF}), so does $F(x)$.
Lemma \ref{BsupGw} leads also to
\begin{equation} \label{dGgeqdw}
d(x) \geq K_1 |x|^{1/2} \, ,
\end{equation}
for some constant $K_1$.

Actually, as soon as the Green function decreases exponentially (with respect to the word distance), these distances are equivalent.
This is the case when $\Gamma$ is non-ameanable (\ref{BsupGwna}) and also for some ameanable groups (for instance the Lamplighter group $\Z \wr \Z_2$, Brofferio \& Woess \cite{brofferiowoess05}).

Remark that (\ref{dGleqdw}) and (\ref{dGgeqdw}) are sharp in the sense that examples exist where the Green function decreases as $\exp (-C|x|)$ or $\exp (-C|x|^{1/2})$ depending on the direction toward infinity (for instance $\Z \wr \Z$ , Revelle \cite{revelle05})

As $B_w(n)\subset B(n)$, the size $V(n)$ of the balls $B(n)$ grows
exponentially, but we prove that the constant in the exponential
(lower and upper) bound is exactly the constant $K$ which appeared
in the definition of the balls.
\begin{pro} \label{voln}
There exist two constants $C_a$ and $C_b$ such that for every integer $n$,
$$
C_a \exp(Kn) \leq V(n) \leq C_b n^3 \exp(Kn) \, .
$$
\end{pro}
\begin{proof}
The proof of the upper bound relies on the following estimate taken from Revelle \cite{revelle05}:
for any finite set $A$,
\begin{equation} \label{sumG}
\sum_{x\in A} G(x) \leq (16 C_e^{-2} +1) (\ln \#A)^3\, .
\end{equation}
By sake of completeness, we give the computation that leads to (\ref{sumG}).
The estimate (\ref{returnP}) (which introduces the constant $C_e$) yields
\begin{eqnarray*}
\sum_{x\in A} G(x) & \leq & \sum_{k=0}^{\infty} \sum_{x \in A} \P [S(k)=x]\\
 & = & \sum_{k\leq C_e^{-3} (\ln \#A)^3 +1}  \sum_{x \in A} \P [S(k)=x] + \sum_{k > C_e^{-3} (\ln \#A)^3 +1}  \sum_{x \in A} \P [S(k)=x]\\
 & \leq &  C_e^{-3} (\ln \#A)^3 +1 +  \# A \sum_{k > C_e^{-3} (\ln \#A)^3 +1} C_e \exp (-C_e k^{1/3}) \\
 & \leq & (16 C_e^{-2} +1) (\ln \#A)^3\, .
\end{eqnarray*}
When $A=B(n)$,  the definition of the ball $B(n)$ and (\ref{sumG})
give
\begin{equation}\label{eqcardBn}
\# B(n) \cdot \exp (-Kn) \leq \sum_{x \in B(n)} G(x) \leq (16
C_e^{-2} +1) (\ln \#B(n))^3\,
\end{equation}
which implies, since $(\ln \#B(n))^3 \leq (\# B(n))^{1/2}$,
that $\# B(n)$ is at most an exponential of $n$. Applying once
more (\ref{eqcardBn}) leads to the $C_b n^3 \exp(Kn)$ upper bound,
for some positive constant $C_b$.

For the lower bound, remark that
$$
\sum_{x\in \pl B(n-1)} G(x) \leq (\# \pl B (n-1)) \, \exp (-K (n-1)) \, .
$$
The left hand side is larger than (or equal to) $1$ since the random walk is transient and it must go through $\pl B(n-1)$.
So with $C_a\fdefeq \exp (-K)$,
$$
V(n) \geq \# \pl B (n-1) \geq C_a \exp (Kn) \, .
$$
\end{proof}

\begin{remk}
For non-ameanable groups, the upper bound in the volume growth
becomes $C_b n \,  \exp(Kn)$. Indeed, using the specific estimate
of the return probability (\ref{returnPna}), the upper bound in
(\ref{sumG}) becomes a constant times $\ln \# A$.
\end{remk}

\section{Limiting shape for the Internal DLA} \label{shapeth}
Let us now state our main result, a comparison between the shape of the Internal DLA and the balls for the hitting distance.
\begin{thm} \label{limshape}
Let $A(\cdot)$ be the Internal DLA on a finitely generated group of exponential growth, associated to a symmetric random walk with finitely supported increment.
Then, for any constants $C_I>3/K$ and $C_O>2$,
$$
\P [ \exists n_0 \mbox{ s.t. } \forall n>n_0 \, , \; B(n-C_I \ln n) \subset A(V(n)) \subset B(n+C_O\sqrt{n}) ] =1 \, .
$$
\end{thm}
This result give an upper bound for the two random variables
$\delta_I(n)$ and $\delta_O(n)$ describing the inner and outer
fluctuation of $A(n)$ with respect to $B(n)$
\begin{eqnarray*}
\delta_I(n) & \fdefeq & n-\inf \{ d(z)/K: z\not\in A(V(n))\}\, ,\\
\delta_O(n) & \fdefeq & \sup \{ d(z)/K: z\in A(V(n))\} -n \, .
\end{eqnarray*}
An alternative definition is that $B(n-\delta_I(n))$ is  the
largest ball  included in $A(V(n))$, and $B(n+\delta_O(n))$ is the
smallest ball which contains $A(V(n))$. Theorem \ref{limshape},
states that for any constants $C_I>3/K$ and $C_O>2$, the events
$\{ \delta_I(n) > C_I \ln n\}$ and $\{ \delta_O(n) >
C_O\sqrt{n}\}$ occur only finitely often, with probability $1$.

We split the proof into two parts, one regarding the inner error
and the other regarding the outer error.

Let us first give some definitions used throughout the proof. Let
$S^j (\cdot )$ be the ${\mathrm j}^{\mathrm{th}}$ random walk that
defines the Internal DLA $A(V(n))$. Recall that all the $S^j
(\cdot)$ start from $e$. Moreover, we let this random walk running
even after it adds a new point to the cluster $A(j-1)$. This
adding time will then be the stopping time:
$$
\sigma^j \fdefeq \inf \{ k\geq 0 :\; S^j(k) \not\in A(j-1) \} \, .
$$
We also define the hitting time of a point $z$ by the (unstopped) random walk $S^j$:
$$
\tau_z^j \fdefeq \inf \{ k\geq 0 :\; S^j(k)=z \} \, .
$$

\subsection{Proof of the inner part}
We will prove that there exists a constant $C_I$ such that, with probability $1$, for all $n$ large enough
$$
B(n - C_I\ln n) \subset A(V(n)) \, .
$$
We start our proof by computing an upper bound for the probability that a given point $z\in B(n)$ does not belong to the Internal DLA at time $V(n)$.
For that purpose, we define the following indicator random variable
$$
N^j=N^j (z) \fdefeq \UN_{\{ \tau_z^j \leq   \sigma^j \}} \, ,
$$
which can be decomposed as $N_j=M_j-L_j$ where
$$
M^j= M^j(z) \fdefeq \UN_{\{ \tau_z^j < \infty \}} \;\; \mbox{and} \;\; L^j=L^j(z) \fdefeq \UN_{\{ \sigma^j \leq \tau_z^j < \infty\}} \, .
$$

Observe that  $z$ does not belong to $A(V(n))$ if and only if
$N^j=0$ for every $j$. So, for any $\lambda
>0$,
\begin{eqnarray*}
\P [z\not\in A(V(n))] & = & \P \left[ \sum_{j=1}^{V(n)} N^j = 0\right] \leq \E \left[ \e^{ -\lambda \sum_{j=1}^{V(n)} N^j} \right]\\
 & = & \E \left[ \e^{-\lambda \sum_{j=1}^{V(n)} M^j + \lambda \sum_{j=1}^{V(n)} L^j} \right]\\
 & \leq & \E \left[ \e^{-2\lambda \sum_{j=1}^{V(n)} M^j} \right]^{1/2} \times \E \left[ \e^{ 2\lambda \sum_{j=1}^{V(n)} L^j} \right]^{1/2}\, .
\end{eqnarray*}
As the $M^j$'s are i.i.d. indicator random variables,
\begin{eqnarray*}
\E \left[ \e^{-2\lambda \sum_{j=1}^{V(n)} M^j} \right] & = & \prod_{j=1}^{V(n)} \E \left[ \e^{-2\lambda  M^j} \right] =  \prod_{j=1}^{V(n)} (1-(1-\e^{-2\lambda}) F(z))\\
 & = & (1-(1-\e^{-2\lambda}) F(z))^{V(n)} \leq \exp (-(1-\e^{-2\lambda}) F(z) V(n))\, .
\end{eqnarray*}
Observe that the $L^j$'s are not independent, but they can be bounded from above by random variables which are independent when knowing where the random walk $S^j$ adds to the cluster.
More precisely, let define the hitting time of $z$ after the adding time $\sigma^j$,
$$
\wt{\tau_z}^j \fdefeq \inf \{k\geq \sigma^j \; \mathrm{s.t.} \; S^j (k)=z\} \, .
$$
Then,
$$
L^j  \leq   \wt{L}^j \fdefeq  \UN_{\{ \wt{\tau_z}^j < \infty \}} \, .
$$
Let $\cG_n$ be the $\sigma$-algebra generated by $S^j (k\wedge \sigma^j)$ for all $j\leq V(n)$ and $k\in \N$, that is the $\sigma$-algebra containing all the information on the $V(n)$ random walks before they add to the cluster.
As the $\wt{L}^j$'s are independent conditioning by $\cG_n$,
\begin{eqnarray*}
\E \left[ \e^{ 2\lambda \sum_{j=1}^{V(n)} L^j} \right] & \leq & \E \left[ \e^{ 2\lambda \sum_{j=1}^{V(n)} \wt{L}^j} \right] = \E \left[ \E \left[ \left. \e^{ 2\lambda \sum_{j=1}^{V(n)} \wt{L}^j} \, \right| \, \cG_n \right] \right]\\
 & = & \E \left[ \prod_{j=1}^{V(n)} \E \left[ \left. \e^{ 2\lambda  \wt{L}^j} \, \right| \, \cG_n \right] \right]  \, .
\end{eqnarray*}
Since the $\wt{L}^j$'s are indicator random variables,
$$
\E \left[ \e^{ 2\lambda \sum_{j=1}^{V(n)} L^j} \right]  \leq   \E \left[ \prod_{j=1}^{V(n)} (1+(\e^{2\lambda}-1) \E [\wt{L}^j \, | \,\cG_n ]) \right]  \leq  \E \left[ \e^{(\e^{2\lambda}-1) \sum_{j=1}^{V(n)} \E [\wt{L}^j \, | \,\cG_n ]} \right] \, .
$$
Finally,
\begin{eqnarray*}
\sum_{j=1}^{V(n)} \E \left[ \left.  \wt{L}^j \, \right| \,\cG_n \right]  & = &  \sum_{j=1}^{V(n)}  \E \left[ \left.  \UN_{\{ \wt{\tau_z}^j < \infty \}} \, \right| \,\cG_n \right] = \sum_{j=1}^{V(n)} \P^{S^j(\sigma^j)} [ \wt{\tau_z}^j < \infty ]\\
 & = &  \sum_{y\in A(V(n))} F(y,z) = \sum_{y\in A(V(n))} F(z,y) \, .
\end{eqnarray*}
The definition of the hitting distance implies that for any $y\in B(z,n)$ and any $y'\not\in B(z,n)$, $F(z,y)>F(z,y')$.
Thus, since $\# A(V(n))=V(n)=\# B(z,n)$,
\begin{equation} \label{optiball}
\sum_{y\in A(V(n))} F(z,y) \leq  \sum_{y\in B(z,n)} F(z,y) = \sum_{y\in B(n)} F(y) \, .
\end{equation}
Putting things together, we get that for any $\lambda >0$, $n$ and $z\in B(n)$,
\begin{equation} \label{inner}
\P [z\not\in A(V(n))] \leq \exp \left( -C_{\lambda} V(n) F(z) +  C_{\lambda}' \sum_{y\in B(n)} F(y) \right) \, .
\end{equation}
with $C_{\lambda} = \frac{1-\e^{-2\lambda}}{2} $ and $C_{\lambda}' = \frac{\e^{2\lambda} -1}{2}$.

We remark that the inequality (\ref{inner}) is valid for every transient groups
whatever its growth rate. For groups of exponential growth,
Proposition \ref{voln}  enables to estimate the last term in the
exponential. Thus, if $z\in B(n-C_I\ln n)$ (that is $F(z)\geq
n^{KC_I} \e^{-Kn}$), using (\ref{sumG}),
$$
\P [z\not\in A(V(n))]  \leq  \e^{-C_{\lambda}
V(n)n^{KC_I}\e^{-Kn}+c \cdot C_{\lambda}'(\ln V(n))^3} \leq
\e^{-C_{\lambda} \cdot C_a n^{KC_I}+ c'\cdot C_{\lambda}'n^3} \leq
\e^{-C_{\lambda}\cdot C_an^3} \, ,
$$
for $C_I>3/K$ and $n$ large enough.
Then,
\begin{eqnarray*}
\sum_{z\in B(n-C_I\ln n)} \P [z\not\in A(V(n))] & \leq & V(n-C_I\ln n) \e^{-C_{\lambda}\cdot C_an^3}\\
 & \leq & C_b (n-C_I\ln n)^3 \e^{Kn-C_I\ln n-C_{\lambda} \cdot C_an^3}
\end{eqnarray*}
which is summable.
Thus the result follows from Borel-Cantelli Lemma.

\subsection{Proof of the outer part}
To get the upper bound for the outer error in Theorem
\ref{limshape}, we will estimate the decay of the expected number
of points in $A(V(n))$ that lies within the spheres $\pl B(n+p)$
for $p\geq 1$. Checking that this decay is exponential in $p^2$
will lead to the desired result.

We define
$$
Z_p(j)\fdefeq  \# (A(j) \cap \pl B(n+p)) =\sum_{i=1}^j \UN_{S^i(\sigma^i) \in \pl B(n+p)} \, ,
$$
and $\nu_p(j)\fdefeq \E[Z_p(j)]$.
If a random walk adds to the cluster in $\pl B(n+p+1)$, that implies it has stayed within the cluster before. It means that
$$
\{ S^i(\sigma^i)\in \pl B(n+p+1)\} \subset \{ \exists x\in \pl B(n+p) \cap A(i-1) \;\; \mathrm{s.t.} \;\; \tau_x^i <\infty\} \, .
$$
Therefore,
\begin{eqnarray*}
\nu_{p+1} (j) & \leq & \sum_{i=1}^j \P [\exists x\in \pl B(n+p) \cap A(i-1) \;\; \mathrm{s.t.} \;\; \tau_x^i <\infty ]\\
 & \leq &  \sum_{i=1}^j \E \left[ \sum_{x\in \pl B(n+p) \cap A(i-1)} \UN_{\tau_x^i <\infty} \right]\\
 & \leq & \sum_{i=1}^j \E \left[ Z_p (i-1) \max_{x\in \pl B(n+p)} \P [\tau_x < \infty] \right] \, ,
\end{eqnarray*}
by independence between $S^{i}$ and $A (i-1)$.
Then,
\begin{equation} \label{recnu}
\nu_{p+1} (j) \leq \sum_{i=1}^j \exp (-K(n+p)) \nu_p (i-1) \, .
\end{equation}
For any fixed $j$, iterating (\ref{recnu}) leads to
$$
\nu_{p} (j) \leq \exp (-Kn(p-1)) \exp (-Kp(p-1)/2) \frac{j^p}{p!} \, ,
$$
by a simple induction on $p$ started from $\nu_1(j)\leq j$ and using
$$
\sum_{i=1}^j \frac{(i-1)^p}{p!} \leq \frac{j^{p+1}}{(p+1)!} \, .
$$
Finally, as $p! \geq p^p \exp(-p)$,
$$
\nu_p(j) \leq  \exp (-Kn(p-1)) \exp (-Kp(p-1)/2) \exp (p) p^{-p} j^p \, .
$$
Thus, for $j=V(n)$ and $p$ larger than some constant independent of $n$,
$$
\nu_p (V(n)) \leq n^{3p} \exp \left( -K \left( \frac{p^2}{3} -n\right) \right) \, .
$$
Therefore, for $C_O > 2$ and $p=C_O \sqrt{n}$,
\begin{eqnarray*}
\P[A(V(n)) \not\subset B(n+C_O\sqrt{n})] & \leq & \P[ Z_{C_O\sqrt{n}}(V(n)) \geq 1] \\
 & \leq & \nu_{C_O\sqrt{n}} (V(n)) \\
 &\leq & \exp \left(- K \left(\frac{C_O^2}{4} -1 \right) n\right) \, ,
\end{eqnarray*}
for $n$ large enough.

Thus, by the Borel-Cantelli Lemma, with probability $1$, for all $n$ large enough,
$$
A(V(n)) \subset B(n+C_O\sqrt{n}) \, .
$$

\section{Lower bounds for the errors when $\cG$ is the homogeneous tree $\T_q$} \label{trees}

In this section we will give lower bounds for the inner and outer
fluctuations in the specific case of Simple Random Walks on the
homogeneous tree.

Let $\cG$ be the homogeneous tree $\T_q$ of degree $q \geq 3$.
When $q$ is even, $\T_q$ is the Cayley graph of the free group
$F_{q/2}$. When $q$ is odd, $\T_q$ is the Cayley graph of the free
product of $F_{(q-1)/2}$ by $\Z_2$. On such a graph, most of the
estimates we need are known exactly (Dynkin \& Maljutov
\cite{dynkin61} and Cartier \cite{cartier72}) for the Simple
Random Walk (nearest neighbours with uniform transition
probability). In particular  since
$$
F(z) =  (q-1)^{-|z|} \, ,
$$
the
hitting distance and the word distance are proportional and the
ball $B_w(n)$ coincides  with $B(n)$ exactly thanks to the use of
the constant $K=\ln (q-1)$. Thus we have
$$
  V(n) = \frac{q(q-1)^n -2}{q-2} \;\mbox{ and }\;\# \pl B(n) =  q (q-1)^n \, .
$$
Let
$$
\xi_n \fdefeq \min \{ k \geq 0 \;\; \mathrm{s.t.} \;\; S
(k) \in \pl B(n)\}\, ,
$$
 be the exit time from the ball $B(n)$. Then, as we deal with Simple Random Walks,
$$
  \P [S(\xi_n)=z] = \frac{1}{q(q-1)^n} \cdot
\delta_{z\in \pl B(n)} \, .
$$

\subsection{Lower bound for the outer error}

\begin{pro} \label{lboutTq}
For any constant $C_o<(\ln (q-1)/\ln q)^{1/2}$,
$$
\P [ A (V(n)) \cap \pl B(n+C_o \sqrt{n} )=\emptyset \; \mbox{infinitely often}]=0\,
$$
\end{pro}

\begin{proof}
By Theorem \ref{limshape}, for any constant $C_{I}$ sufficiently large, if one defines the increasing sequence of events
$$
\Omega_{n_{0}} \fdefeq \left\{ \forall n\geq n_{0}, \,
B(n-C_{I}\ln(n)+1)  \subset A(V(n))\right\} \, .
$$
One has
\begin{equation} \label{innerr}
\lim_{n_0 \rightarrow \infty} \P [\Omega_{n_{0}}]=1 \, .
\end{equation}
Suppose that $\Omega_{n_{0}}$
occurs and fix $n> n_{0}$. Denote $R\fdefeq n-1-C_I \ln (n-1)$.
Let the first $V(n-1)$ random walks
build $A(V(n-1))$, which covers $B(R)$. Then we consider the
next $V(n)-V(n-1)$ random walks and stop them as soon as they
reach $\partial B(R)$. Let $N_{x}$ be the number of random walks
that are stopped in $x\in\partial B(R)$. Letting the stopped
random walks start again we can reconstruct the original Internal
DLA. The same random walks construct new smaller Internal DLA's,
$A_{x}(N_{x})$, that start in every $x\in \pl B(R)$ and are built
up using $N_{x}$ random walks . We state that if $B(R)\subset
A(V(n-1))$ then
$$
A(V(n))\supset A(V(n-1)) \bigcup_{x\in\partial B(R)}A_{x}(N_{x}) \, .
$$
Thus, if there exists $x\in\partial B(R)$ such that
$A_{x}(N_{x})\cap\partial B(R+r)\neq \emptyset$, then
$A(V(n))\cap\partial B(R+r)\neq \emptyset$. Observe also that the
random walks that start on $\partial B(R)$ are independent and,
thus, all the $\left\{ A_{x}(N_{x})\right\} _{x\in\partial B(R)}$
are independent once one knows all the $N_{x}$, that is  they are
independent for the probability $\P [ \, \cdot \; | \; \sigma_R]$.
Here $\sigma_R$ is the $\sigma$-algebra generated by the
$S^i(k\wedge \xi_R)$ for all $i\leq N\fdefeq V(n)-V(n-1)$ and all
$k\in \N$, that is the $\sigma$-algebra containing all the
information on the $N$ random walks until they hit $\pl B(R)$.

For all $x\in\partial B(R),$ there exits a path $\cC_{x}$ of (graph) length
$r$ from $x$ to $\pl B(R+r)$.
Note that
$$
\{ A_{x}(N_{x})\cap\partial B(R+r)=\emptyset \} \subset \{ \cC_{x}\not\subset A_{x}(N_{x})\} \, .
$$
We split the $N_{x}$ random walks that start at $x$ into $\left\lfloor \frac{N_{x}}{r}\right\rfloor $
packets of size $r$ ($\left\lfloor \cdot \right\rfloor$ is the lower integer part), the rest being ignored.
The probability that all the random walks of one packet follow exactly $\cC_x$ when they start, is $(q^{-r})^r$.
By definition of the Internal DLA, this event guaranty that $\cC_x$ is contained in $A_x(N_x)$.
Thus, since all packets are independent
$$
\P [\cC_{x}\not\subset A_{x}(N_{x}) \; | \; \sigma_R]
\leq(1-q^{-r^{2}})^{\left\lfloor \frac{N_{x}}{r}\right\rfloor }\leq(1-q^{-r^{2}})^{\frac{N_{x}}{r}-1} \, .
$$
Finally, for  $n>n_0$,
\begin{eqnarray*}
\P [\left\{ A(n)\cap\partial B(R+r)= \emptyset  \right\} & \cap & \Omega_{n_{0}}]\\
 & \leq & \P [\forall x\in\partial B(R), \, A_{x}(N_{x})\cap\partial B(R+r)=\emptyset ]\\
 & = & \E \left[ \prod_{x\in\partial B(R)}\P [A_{x}(N_{x})\cap\partial B(R+r)= \emptyset \; | \; \sigma_R]  \right]\\
 & \leq & \E \left[ \prod_{x\in\partial B(R)}(1-q^{-r^{2}})^{\frac{N_{x}}{r}-1}\right]\\
 & = & \E \left[ (1-q^{-r^{2}})^{\sum_{x\in\partial B(R)}\frac{N_{x}}{r}-1} \right]\\
 & = & (1-q^{-r^{2}})^{\frac{V(n)-V(n-1)}{r}-\# \partial B(R)}\\
 & \leq & \exp \left( -q^{-r^{2}}\left(\frac{V(n)-V(n-1) }{r} -\#\partial B(R)\right)\right) \, .
\end{eqnarray*}
Remark that taking $r=C_o \sqrt{n}+1+C_I \ln (n-1)$, with a
constant $C_o<(\ln (q-1)/ \ln q)^{1/2}$ and $C_I >\frac{1}{2}$, we
get
$$
\exp \left( -q^{-r^{2}}\left(\frac{V(n)-V(n-1) }{r} -\#\partial
B(R)\right) \right) \leq \exp (-\exp (cn)) \, ,
$$
for some constant $c>0$.
Then, by Borel-Cantelli Lemma,
$$
\P [\{A(V(n)) \cap \pl B(n+C_o \sqrt{n}) \; \mathrm{i.o.} \} \cap \Omega_{n_0} ] =0 \, .
$$
Letting $n_0$ goes to infinity leads (by (\ref{innerr})) to the result.
\end{proof}

\subsection{Lower bound for the inner error}

For the inner error, we prove a weaker result than the one on the
outer error. Namely, we get that infinitely often, the inner error
is of order $\ln n$. It assures the sharpness of the order $\ln n$
for a general upper bound but does not prevent the inner error to
be sometimes smaller. Thus our result  is somehow weaker than the
one announced by Etheridge and Lawler.
\begin{pro} \label{lbinnTq}
For any constant $C_i<1/(2\ln (q-1))$,
$$
\P [A(V(n))^c \cap B(n-C_i \ln n) \neq \emptyset \; \mbox{infinitely often}]=1\, .
$$
\end{pro}

\begin{proof}
Let $N \fdefeq V(n)$ and $R \fdefeq n-1-C_i \ln n$.
Then the event $\{A(N)^c \cap B(R+1) \neq \emptyset \}$ contains the event that there exists a point on the sphere of radius $R$ attained by none of $N$ independent random walks.

Let fix $r<R$, for every $z\in \pl B(R)$, we define its class as
$$
[z] \fdefeq \pl B(R) \cap B(z,2r) \, ,
$$
and let $A_R \fdefeq \bigcup_{i=1}^N [S^i (\xi_R)]$ be the set of
all classes that are hit by the $N$ random walks when they first
exit $B(R)$. Then we have
\begin{eqnarray*}
\P[A(N)^c & \cap & B(R) \neq \emptyset] \\
 & \geq &\P [\exists z \in \pl B(R) \;\; \mathrm{s.t.} \;\;  \forall i\leq N , \, \tau_z^i =\infty ]\\
 & \geq & \P [ \exists z \in \pl B(R) \;\; \mathrm{s.t.} \;\; \forall i\leq N , \, \tau_z^i =\infty \;\; \mathrm{and} \;\; A_R \varsubsetneq \pl B(R)]] \\
 & = & \!\!\!\! \sum_{A \varsubsetneq \pl B(R)} \!\! \P [ \exists z \in \pl B(R) \; \mathrm{s.t.} \; \forall i\leq N , \, \tau_z^i =\infty \; | \; A_R=A] \P [A_R=A] \, .
\end{eqnarray*}
To get a lower bound for the above conditional probability we will
need, once the set $A$ is fixed, to choose a suitable point called
$z_A$ on $\pl B(R)$ that has few chances to be hit by the random walks
after $\xi^i_R$. For that purpose, we construct the following
algorithm. Suppose $A \varsubsetneq \pl B(R)$ is a fixed set, then
we will write a ``Mouse algorithm'' whose aim is to find the most
secure place $z_A$ on $\pl B(R)$ for a mouse knowing that $\# A$
cats will start simple random walks (hunts) from each of the
points in $A$.

The mouse starts from the origin and proceeds toward $\pl B(R)$
choosing the direction of the subtree with the fewest cats. Denote
$(m_j)$ ($j=0$ to $R$) the successive positions of the mouse, that
is the steps of our algorithm:
\begin{itemize}
\item[$\bullet$] $m_0=e$;
\item[$\bullet$] Choose $m_{j+1}$ inside $\{y\sim m_j \, : \; |y|=j+1$ and $B(y,R-j)\cap A$ is minimal$\}$.
\end{itemize}
Let $z_A \fdefeq m_{R+1}$.
It is easy to check that $z_A \not\in A$.
Recall that $\sigma_R$ is the $\sigma$-algebra generated by the $S^i(k\wedge \xi_R)$ for all $i\leq N$ and all $k\in \N$, that is the $\sigma$-algebra containing all the information on the $N$ random walks before they hit $\pl B(R)$.
Then,
\begin{eqnarray*}
\P [\exists z \in \pl B(R) & \mathrm{s.t.} & \forall i\leq N , \, \tau_z^i=\infty \; | \; A_R=A]\\
 & \geq & \P [ \forall i\leq N , \, \tau_{z_A}^i=\infty \; | \; A_R=A]\\
 & = & \P [ \forall i\leq N , \, \xi_R^i \leq \tau_{z_A}^i=\infty \; | \; A_R=A]\\
 & = & \E \left[ \P [ \forall i\leq N , \, \xi_R^i \leq \tau_{z_A}^i=\infty \; | \; \sigma_R ]\; | \; A_R=A \right] \\
 & = & \E \left[ \left. \prod_{i=1}^N \P [ \xi_R^i \leq \tau_{z_A}^i=\infty \; | \; \sigma_R ]\; \right| \; A_R=A \right] \\
 & = & \E \left[ \left. \prod_{i=1}^N \P^{S^i (\xi_R)} [ \tau_{z_A}=\infty ]\; \right| \; A_R=A \right] \, .
\end{eqnarray*}
By transience of the random walk, $\P^{S^i (\xi_R)} [ \tau_{z_A}=\infty ]=1- F(S^i (\xi_R) , z_A)$ has a uniform (in all parameter) strickly positive lower bound.
Thus, there exists a constant $c>0$ such that
\begin{eqnarray*}
\E \left[ \left. \prod_{i=1}^N \P^{S^i (\xi_R)} [ \tau_{z_A}=\infty ]\, \right| \, A_R=A \right] & \geq & \E \left[ \left.\exp \left( -c \sum_{i=1}^N F(S^i (\xi_R) , z_A)\right) \, \right| \, A_R=A \right] \\
 & \geq & \exp \left( -c \sum_{i=1}^N \E \left[ F(S^i (\xi_R) , z_A) \; | \; A_R=A \right] \right) \\
 & = & \exp \left( -c N \E \left[ F(S^1 (\xi_R) , z_A) \; | \; A_R=A \right] \right)\, .
\end{eqnarray*}
By definition of the Mouse algorithm, we easily check that, for every $j$,
$$
T_j \fdefeq B(z_A, 2(R+1-j))\cap A = B(m_j , R+1-j) \cap A \, ,
$$
and that its cardinality is $\# T_j \leq \# A (q-1)^{-j} $. Remark
also that since $A$ must be a union of sets of the type $\pl B(R)
\cap B(z,2r)$ with $z\in \pl B(R)$ (otherwise $\P [A_R=A]=0$),
$z_A\not\in A$ implies  $d(z_A,A) > 2r$. Then,
\begin{eqnarray*}
\E \left[ F(S^1 (\xi_R) , z_A) \; | \, A_R=A\right] & = &
 \!\! \sum_{j=0}^{R-r} \E \left[ \UN_{S^1(\xi_R) \in T_j \bk T_{j+1}} F(S^1 (\xi_R) , z_A) \; | \, A_R=A\right]\\
 & = &  \!\! \sum_{j=0}^{R-r} \E \left[ \UN_{S^1(\xi_R) \in T_j \bk T_{j+1}} (q-1)^{-2(R+1-j)} \; | \, A_R=A\right]\\
 & = &  \!\! \sum_{j=0}^{R-r}  (q-1)^{-2(R+1-j)} \P \left[ S^1(\xi_R) \in T_j \bk T_{j+1}  \; | \, A_R=A\right] \, .
\end{eqnarray*}

Remark that the event $A_R =A$ is invariant by any bijection $\gamma$ of $\pl B(R)$ such that $\gamma (A)=A$ and for any $x\in \pl B(R)$, $[\gamma (x)]=\gamma ([x])$.
Thus, the law of $S^1 (\xi_R)$ under the condition $A_R=A$ is the uniform distribution on $A$.
Thus,
$$
\P \left[ S^1(\xi_R) \in T_j \bk T_{j+1}  \; | \; A_R=A\right] = \frac{\# (T_j \bk T_{j+1})}{\#A} \leq \frac{\# T_j}{\#A} \, .
$$
Therefore,
$$
\E \left[ F(S^1 (\xi_R) , z_A) \; | \; A_R=A\right] \leq \sum_{j=0}^{R-r}  \frac{(q-1)^{-2(R+1-j)}}{(q-1)^j} \leq (q-1)^{-R-r} \, .
$$
Finally,
\begin{eqnarray*}
\P [\exists z\in \pl B(R) \;\; \mathrm{s.t.} \;\; \forall i \leq N ,\, \tau_z^i=\infty ] & \geq & \sum_{A \varsubsetneq \pl B(R)} \! \exp \left( -cN (q-1)^{-R-r} \right) \P [A_R=A] \\
 & = & \exp \left( -cN (q-1)^{-R-r} \right) \P [A_R \varsubsetneq \pl B(R)] \, .
\end{eqnarray*}
Remark that $\P [A_R \varsubsetneq \pl B(R)]$ is the probability
that $N$ balls uniformly distributed into $M=q(q-1)^{R-r}$ boxes
leave at least one empty box. By Lemme \ref{multi}, that will be
prove at the end of the section,
$$
\P [A_R \varsubsetneq \pl B(R)] \geq 1- \exp \left( -\frac{q(q-1)^{R-r}}{2} \exp \left( -4N \frac{(q-1)^{r-R}}{q} \right) \right) \, .
$$
When $N=V(n)\leq q(q-1)^n$, $R=n-C_i \ln n$ and $r=C_1 \ln n$ (the constant $C_1$ will be fixed below), we easily check that
\begin{multline*}
\exp \left( -\frac{q(q-1)^{R-r}}{2} \exp \left( -4N \frac{(q-1)^{r-R}}{q} \right) \right) \\
\leq \exp \left( -\frac{q\, n^{-(C_1+C_i)\ln (q-1)}}{2} \exp \left( n\ln (q-1) -4n^{(C_1+C_i)\ln (q-1)} \right) \right) \, .
\end{multline*}
So, with $C_1=\frac{1}{2\ln (q-1)}$, for any constant $C_i<\frac{1}{2\ln (q-1)}$,
$$
\P [A_R \varsubsetneq \pl B(R)] \build{\longrightarrow}_{n \rightarrow \infty}^{} 1 \,  .
$$
With the same constants $C_1$ and $C_i$, as $C_1>C_i$,
$$
\exp \left( -cN (q-1)^{-R-r} \right)
\geq \exp \left( -c\cdot q \cdot n^{-(C_1-C_i)\ln (q-1)} \right) \, ,
$$
also tends to $1$ when $n$ goes to infinity.
Finally,
$$
 \P [\exists
z\in \pl B(n-C_i \ln n) \;\; \mathrm{s.t.} \;\; \forall i \leq
V(n) ,\, \tau_z^i=\infty ] \build{\longrightarrow}_{n \rightarrow
\infty}^{} 1 \, .
$$
Thus $\P [A(V(n))^c \cap B(n-C_i \ln n) \neq \emptyset
]\rightarrow 1$ which completes the proof by Fatou Lemma.
\end{proof}

We now prove the estimate on the multinomial law that we used in
the proof.
\begin{lem} \label{multi}
Let put independently $N$ balls into $M$ boxes ($M\geq 3$) with
uniform probability. Let $\P_{N,M}$ denote the probability law of
that process. Then,
$$
\P_{N,M} [\mbox{Each box contains at least one ball}] \leq \exp
\left( -\frac{M}{2} \exp \left( -4 \frac{N}{M} \right) \right) \,
.
$$
\end{lem}

\begin{proof}
Let $(b_k)$ ($k=1$ to $M$) denote the $M$ boxes and $(a_i)$ ($i=1$
to $N$) the $N$ balls. Let also define
\begin{eqnarray*}
p_{N,M} & \fdefeq & \P_{N,M} [ \forall k \leq M , \, \exists i\leq N \;\; \mathrm{s.t.} \;\; a_i \in b_k ] \; ;\\
q_{N,M} & \fdefeq & \P_{N,M} [\exists i \leq N \;\; \mathrm{s.t.} \;\; a_i \in b_1] \; ; \\
I_1 & \fdefeq & \{ i \;\; \mathrm{s.t.} \;\; a_i \in b_1 \} \, .
\end{eqnarray*}
Therefore,
\begin{eqnarray*}
p_{N,M} & = & \sum_{I\neq \emptyset} \P_{N,M} [\{\forall k\neq 1 \, \exists i \not\in I \;\; \mathrm{s.t.} \;\; a_i \in b_k \} \cap I_1=I]\\
 & = & \sum_{I\neq \emptyset} \P_{N,M} [\forall k\neq 1 \, \exists i \not\in I \;\; \mathrm{s.t.} \;\; a_i \in b_k \; | \;  I_1=I]\times \P_{N,M} [I_1=I]\\
 & = & \sum_{I\neq \emptyset} p_{N-\# I,M-1} \times  \P_{N,M} [I_1=I]\\
 & \leq  & \sum_{I\neq \emptyset} p_{N,M-1}  \times \P_{N,M} [I_1=I]\\
 & = & p_{N,M-1} \cdot q_{N,M} \, .
\end{eqnarray*}
The third equality comes from the fact that conditioning by
$I_1=I$, the $N-\#I$ remaining balls are uniformly distributed
among the $M-1$ remaining boxes. The inequality comes from the
monotonicity in $N$ of $p_{N,M}$. By induction,
$$
p_{N,M} \leq \prod_{j=2}^M q_{N,j}\cdot p_{N,1} \, .
$$
As $p_{N,1}=1$ and $q_{N,j} =1 - \left( 1-\frac{1}{j} \right)^N$,
\begin{eqnarray*}
p_{N,M} & \leq & \prod_{j=\lfloor (M+1)/2 \rfloor}^M q_{N,j} \leq \prod_{j=\lfloor (M+1)/2 \rfloor }^M \left( 1- \left( 1-\frac{2}{M} \right)^N \right)\\
  & \leq & \left(1- \left( 1-\frac{2}{M} \right)^N \right)^{M/2} \leq  \left( 1-\exp \left( -c \frac{2N}{M} \right) \right)^{M/2}\\
 & \leq & \exp \left( -\frac{M}{2} \exp \left( -c\frac{2N}{M} \right) \right)\, ,
\end{eqnarray*}
for some constant $c>0$ ($c=2$ for instance, since $M\geq 3$).
\end{proof}

\bibliographystyle{amsplain}
\providecommand{\bysame}{\leavevmode\hbox to3em{\hrulefill}\thinspace}
\providecommand{\MR}{\relax\ifhmode\unskip\space\fi MR }
\providecommand{\MRhref}[2]{%
  \href{http://www.ams.org/mathscinet-getitem?mr=#1}{#2}
}
\providecommand{\href}[2]{#2}

\begin{flushleft}
S\'ebastien Blach\`ere\\
Universit\'e Aix-Marseille 1\\
CMI, 39 rue Joliot-Curie\\
13453 Marseille Cedex, France
\end{flushleft}
\begin{flushleft}
Sara Brofferio\\
Universit\'e Paris-Sud\\
Laboratoire de Math\'ematiques, b\^at. 425,\\
91405 Orsay Cedex, France.
\end{flushleft}

\end{document}